# An Inductive Proof of Bertrand's Postulate

Bijoy Rahman Arif
University of Nebraska, Omaha, NE
bijoyarif71@yahoo.com

## Abstract

In this paper, we are going to prove a famous problem concerning prime numbers. Bertrand's postulate states that there is always a prime $p$ with $n < p < 2n$, if $n > 1$. Bertrand's postulate isn't a newer one to be proven, in fact, after Bertrand's assumption and numerical evidence, Chebyshev was the first person who proved it. Subsequently, Ramanujan proved it using properties of Gamma function, $\Gamma(x)$ , and Erdös published a simpler proof with the help of Primorial function, *p#*. Our approach is unique in the sense that we have used mathematical induction for finding the upper and lower bounds for the second Chebyshev function, $\psi(x)$ and they are even stronger than Ramanujan's bounds finding using Gamma function, $\Gamma(x)$ . Otherwise, our approach is similar the way Ramanujan proved it.

## Definition

We define function $\upsilon(x)$ and $\psi(x)$ conventionally [1] as:

$$\upsilon(x) = \sum_{p \le x} \log(p) = \log \prod_{p \le x} p$$

$$\psi(x) = \sum_{p^m \le x} \log(p)$$

Since $p^2 \le x, p^3 \le x, \ldots$ are equivalent to $p \le x^{1/2}, p \le x^{1/3}, \ldots$ , we have [1] , [2]:

$$\psi(x) = \upsilon(x) + \upsilon(x^{1/2}) + \upsilon(x^{1/3}) + \ldots = \sum_{m \ge 1} \upsilon(x^{1/m})$$

and so $\psi(2n) = \upsilon(2n) + \upsilon((2n)^{1/2}) + \upsilon((2n)^{1/3}) + \ldots = \sum_{m \ge 1} \upsilon((2n)^{1/m})$

## Proof

We know [2], $\log((2n)!) = \psi(2n) + \psi(n) + \psi(\frac{2n}{3}) + \ldots$ , ………………. (1)

From relation of $\upsilon(x)$ and $\psi(x)$ , we get [2] :

$$\psi(2n) - 2\psi(\sqrt{2n}) = \upsilon(2n) - \upsilon((2n)^{1/2}) + \upsilon((2n)^{1/3}) - \ldots \quad , \text{………………. (2)}$$

Let $N_n = \frac{(2n)!}{n!n!}$ , then from (1) :

$$\log(N_n) = \psi(2n) - \psi(n) + \psi(\frac{2n}{3}) - \ldots \quad , \text{………………. (3)}$$

As $\upsilon(x)$ and $\psi(x)$ are steadily increasing function, we find from (2) and (3) that:

$$\psi(2n) - 2\psi(\sqrt{2n}) \le \upsilon(2n) \le \psi(2n) \quad , \text{………………. (4)}$$

and $\psi(2n)-\psi(n) \leq \log(N_n) \leq \psi(2n)-\psi(n)+\psi(\frac{2n}{3})$ , ………………. (5)

Now $N_{n+1} = \frac{(2n+2)!}{(n+1)!(n+1)!} = 2.\ \frac{2n+1}{n+1}\ .N_n$

for $n = 1$, $\frac{2n+1}{n+1} = \frac{3}{2}$ and $\lim_{n\to\infty}\frac{(2n+1)}{(n+1)} = 2$ which implies for $n \geq 1$,

$$3N_n \leq N_{n+1} \leq 4N_n \text{ , ………………. (6)}$$

We assume, $a=\log 3$ and $b=\log 4$

for $n = 2$, $N_2 = \frac{2!}{1!1!} = 2$; $\log 2 < b\times 2$

and $n = 5$, $N_5 = \frac{10!}{5!5!} = 252$; $a\times 5 < \log 252$

We assume, $\log(N_n) < bn$ if $n \geq 2$, and $an < \log(N_n)$ if $n \geq 5$, ………………. (7)

It follows from (6) and (7) that: $\log 3 + \log(N_n) \leq \log(N_{n+1}) \leq \log 4 + \log(N_n)$ implies $a(n+1) < \log(N_{n+1}) < b(n+1)$ , by induction (7) is proven.

It follows from (5) and (7) that:

$\psi(2n)-\psi(n) < bn$ if $n \geq 2$, ………………. (8)

$\psi(2n)-\psi(n)+\psi(\frac{2n}{3}) > an$ if $n \geq 5$, ………………. (9)

Now, changing $n$ to $\frac{n}{2},\frac{n}{4},\frac{n}{8},\ldots$ in (8) and adding up all the results, we get:

$\psi(2n) < 2bn$ if $n \geq 2$, ………………. (10)

Finally we have from (4) and (10):

$$\psi(2n)-\psi(n)+\psi(\frac{2n}{3}) < \upsilon(2n)+2\psi(\sqrt{2n})-\upsilon(n)+\psi(\frac{2n}{3})$$

$< \upsilon(2n)-\upsilon(n)+2b\sqrt{2n}+\frac{2b}{3}n$ , ………………. (11)

We conclude from (9) and (11) that:

$\upsilon(2n)-\upsilon(n) > an-2b\sqrt{2n}-\frac{2b}{3}n > 0$ if $n > 505$, by considering right hand side a quadratic equation.

It is easy to verify using simple computer program there are primes for $1<n\leq 505$ , hence, we have proved Bertrand's postulate.